\documentclass[12pt]{article}
\usepackage{amssymb}

\newtheorem{theorem}{Theorem}
\newtheorem{lemma}{Lemma}

\newcommand{\eps}{\varepsilon}
\renewcommand{\phi}{\varphi}

\title{On a certain generalization of the Balog-Szemeredi-Gowers theorem}
\author{Ernie Croot and Evan Borenstein}

\begin{document}

\maketitle

The Balog-Szemer\'edi-Gowers theorem has a rich history, and is one
of the most useful tools in additive combinatorics.  It began with the
a paper by Balog and Szemer\'edi \cite{balog}, and then was refined
by Gowers \cite{gowers} to the following basic result (actually, 
Gowers proved somewhat more than we bother to state here):

\begin{theorem} \label{BSG}  There exists an absolute constant $\kappa > 0$
such that the following holds for all finite subsets $X$ and $Y$ 
of size $n > n_0$ of an abelian group:  Suppose that there are at least 
$C n^3$ solutions to $x_1 + y_1 = x_2 + y_2$, 
$x_i \in X$ and $y_i \in Y$.  Then, $X$ contains a subset $X'$, 
of size at least $C^\kappa n$, such that
$$
|X' + X'|\ \leq\ C^{-\kappa} n.
$$
\end{theorem}

Sudakov, Szemer\'edi and Vu \cite{sudakov} proved a refinement of
this theorem (Balog \cite{balog2} independently obtained a  
similar result), given as follows:

\begin{theorem} \label{SSV} 
Let $n,C,K$ be positive numbers, and let $A$ and $B$
be two sets of $n$ integers.  Suppose that there is a bipartite
graph $G(A,B,E)$ with at least $n^2/K$ edges and 
$|A +_G B| \leq Cn$.  Then one can find a subset $A' \subset A$ and a 
subset $B' \subset B$ such that $|A'| \geq n/16K^2$, 
$|B'| \geq n/4K$ and $|A'+B'| \leq 2^{12} C^3 K^5 n$.  
\end{theorem}

\noindent {\bf Remark.}  It is not difficult to show that this theorem,
along with some lemmas and theorems of Ruzsa (the Ruzsa triangle 
inequality \cite{tao}, and the Ruzsa-Plunnecke Theorem 
\cite{ruzsa}), implies that we may take $\kappa < 20$ in Theorem
\ref{BSG}. 
\bigskip

In the same paper, Sudakov, Szemer\'edi and Vu \cite[Theorem 4.3]{sudakov}
proved the following powerful hypergraph version of the 
Balog-Szemer\'edi-Gowers Theorem: 

\begin{theorem} For any positive integer $k$, there are polynomials
$f_k(x,y)$ and $g_k(x,y)$ with degrees and coefficients depending only
on $k$, such that the following holds.  Let $n,C,K$ be positive numbers.
If $A_1, ..., A_k$ are sets of $n$ positive integers, $H(A_1,...,A_k,E)$
is the $k$-partite, $k$-uniform hypergraph with at least $n^k/K$
edges, and $|\oplus_{H i =1}^k A_i| \leq Cn$, then one can find subsets 
$A_i' \subset A_i$ such that 

$\bullet$ $|A_i'| \geq n/f_k(C,K)$ for all $1 \leq i \leq k$;

$\bullet$ $|A_1' + \cdots + A_k'| \leq g_k(C,K) n$.
\end{theorem}
The notation $\oplus_H$ means that the sum is restricted to the hypergraph
$H$.
\bigskip

Beautiful and useful as it is, it would be nice if one had some 
control on the degrees of these polynomials $f$ and $g$.
And, for particular applications that we 
(Croot and Borenstein) have in mind, it would be good to be able to 
control the rate of growth of sums $A_1' + \cdots + A_\ell'$, where
$\ell$ is much smaller than $k$ -- it would be good to 
be able to bound the size of this sum from above by 
\begin{equation} \label{implicit}
C^{1 + \eps} K^{d_k} n,
\end{equation}
where $d_k$ depends only on $k$.  Perhaps such a bound can be 
developed by modifying the proof of Sudakov, Szemer\'edi and Vu; 
however, in the present paper, we take a different tack, and produce 
an alternate proof of a related hypergraph Balog-Szmeredi-Gowers 
theorem, where such an upper bound as (\ref{implicit}) will be implicit,
though only for the case where $A_1 = \cdots = A_k$.  
In our proof, we will use some of the same standard tricks 
as Sudakov, Szemer\'edi and Vu do in their proof.

The notation we use to describe this theorem, and its proof, will be
somewhat different from that used by Sudakov, Szemer\'edi and Vu.
Furthermore, we will not attempt here to give the most general formulation
of the theorem.

\begin{theorem} \label{main_theorem}  
For every $0 < \eps < 1/2$ and $c > 1$, there exists $\delta > 0$,
such that the following holds for all $k$ sufficiently large, and all
sufficiently large finite subsets $A$ of an additive abelian group:  
Suppose that 
$$
S\ \subseteq\ A \times A \times \cdots \times A\ =\ A^k,
$$
and let 
$$
\Sigma(S)\ :=\ \{ a_1 + \cdots + a_k\ :\ (a_1,...,a_k) \in S\}.
$$
If 
$$
|S|\ \geq\ |A|^{k-\delta},\ {\rm and\ } |\Sigma(S)|\ <\ |A|^c,
$$
then there exists 
$$
A'\ \subseteq\ A,\ |A'|\ \geq\ |A|^{1-\eps},
$$
such that
$$
|\ell A'|\ =\ |A' + \cdots + A'|\ \leq\ |A'|^{c(1 + \eps\ell)}.
$$ 
\end{theorem}

\section{Proof of Theorem \ref{main_theorem}}

\subsection{Notation and basic assumptions}

It will be advantageous to describe the proof in terms of strings.
So, the set 
$$
S\ \subseteq\ A^k
$$
will be thought of as a collection of strings of length $k$:
$$
x_1 x_2 \cdots x_k,
$$
where each $x_i \in A$. 

Often, we split these strings up into substrings; for example, 
the string
$$
x\ =\ x_1 \cdots x_k
$$
can be written as a product of a ``left substring $\ell$ of length $k/2$''
(assume $k$ is even) and a ``right substring $r$ of length $k/2$''.  So,
$$
x\ =\ \ell r.
$$
\bigskip

We may assume that 
$$
k\ =\ 2^n,
$$
since if this is not the case, then we let $k'$ be the largest power
of $2$ of size at most $k$, and proceed as follows:  Given a string
$x_1\cdots x_k$ in $S$, we write it as a product $\ell_x r_x$, where 
$$
\ell_x\ :=\ x_1 \cdots x_{k'}\ \ {\rm and\ \ } r_x\ :=\ x_{k'+1} \cdots x_k. 
$$
Now, for some string $y$ we will have that $r_x = y$ for at   
least $|S|/|A|^{k-k'}$ choices for $x \in S$.  Letting $S'$ denote
the set of all strings $\ell_x$ with $r_x = y$, we will have
$$
|S'|\ \geq\ |A|^{k'-\delta},
$$
and clearly
$$
|\Sigma(S')|\ \leq\ |\Sigma( \{ \ell_x y\ :\ x \in S'\})|\ 
\leq\ |\Sigma(S)|\ \leq\ |A|^c.
$$
 
So, we could just assume that our $k$ had this value $k'$ all along
(remember, we get to choose $k$ to be as large as needed to get the
desired conclusion).

\subsection{The suppression of subscripts, and a comment about iteration}

In the proof of our theorem, we will iteratively replace our initial
set $S$ with other, smaller and smaller sets having certain useful 
properties.  If we were so inclined, we could describe this iteration
by saying that we produce a sequence of sets
$$
S_0 := S,\ S_1,\ S_2, ...,\ S_t,\ {\rm where\ } S_i\ \subseteq\ A^{k_i},
\ |S_i|\ \geq\ |A|^{k_i - \delta_i}.
$$
The trouble with this is that it leads to a proliferation of subscripts,
which can be unpleasant.  

Instead of introducing subscripts, we use the ``assignment operator'',
denoted by
$$
S\ \leftarrow\ S',
$$
which means that the set $S$ gets ``reassigned'' to the set $S'$.  
So, it is worth keeping in mind that later into the proof, 
$S$ refers to a different set than at the start of the proof.  
The same will be true of $k$ and $\delta$.

\subsection{Lengths of iterations and the choice of $\delta$ and
$k$}

At almost every step of our iteration, we will replace $S \subseteq A^k$ with
$S' \subseteq A^{k'}$, satisfying  
$$
|S'|\ \geq\ A^{k' - \delta},\ {\rm and\ }
|A|^{1-O(\delta)}\ \leq\ |\Sigma(S')|\ \leq\ |\Sigma(S)|^{1- \eps/400c}
$$
Clearly, for $\delta > 0$ small enough, the number of such iterations
we can take will be bounded from above in terms of $\eps$ and $c$.
Furthermore, since at each step, $k'$ is at least half the size of
$k$, so long as the initial value of $k$ is large enough in terms 
of $c$ and $\eps$, we will not run out of dimensions. 

Since our theorem is a qualitative result, in that it does not 
even attempt to explain how $\delta$ or $k$ 
depends on $\eps$ and $c$, there 
is no need to be more precise about just how small one needs take
$\delta$ or how large to take $k$, 
in order for our iteration process to terminate and prove our theorem.   

\subsection{The iteration part of the argument} \label{iteration_section}

Given a string $x$ of length $k/2$, we let $R_x$ denote the set of 
all strings $y$ of length $k/2$ such that 
$$
xy\ \in\ S.
$$
We analogously define $L_x$.
\bigskip

We will now select an $x$, and therefore $R_x$, very carefully, so 
that it satisfies certain useful properties:  We begin with the
inequality
$$
\sum_x |R_x|\ =\ |S|\ \geq\ |A|^{k - \delta}.
$$  
We now apply the following lemma, which is easily proved upon using
the Cauchy-Schwarz inequality:

\begin{lemma} \label{intersect_lemma}  Suppose that $V$ is a set of 
$n$ elements, and suppose that 
$$
U_1, U_2, ..., U_r\ \subseteq\ V
$$
satisfy
$$
\sum_{i=1}^r |U_i|\ \geq\ r n^{1-\delta}.
$$
Then, there exists $1 \leq j \leq r$ such that 
$$
\sum_{1 \leq i \leq r} |U_i \cap U_j|\ \geq\ r n^{1-2\delta}.
$$
\end{lemma}
\bigskip

From this lemma we easily deduce that there exists $x$ such that  
$$
\sum_y |R_x \cap R_y|\ \geq\ |A|^{k - 2\delta}.    
$$
Next, we let 
\begin{equation} \label{S'U}
S'\ :=\ \{ yz \in S\ :\ z \in R_x\},
\end{equation}
and we observe that 
$$
|S'|\ =\ \sum_y |R_x \cap R_y|\ \geq\ |A|^{k - 2\delta};
$$
so, $S'$ is not too much smaller than $S$.  

We now make a reassignment:
$$
S\ \leftarrow\ S',\ \delta\ \leftarrow\ 2 \delta,
$$
and observe that $S$ now satisfies
$$
|S|\ \geq\ |A|^{k - \delta},
$$
and we in addition have that every element of $S$ can be expressed as 
$yz$, where $z \in R_x$.
\bigskip

Now suppose that there is a string $y$ of length $k/2$ such that if
$$
|R_y|\ \geq\ |A|^{k/2 - 2 \delta},
$$
then
$$
|\Sigma(R_y)|\ \leq\ |\Sigma(S)|^{1-\eps/400c}. 
$$
If this occurs, then we make another reassignment: 
$$
S\ \leftarrow\ R_y,\ k\ \leftarrow\ k/2,\ \delta\ \leftarrow\ 2\delta,
$$ 
and we start back at the very beginning of this subsection 
\ref{iteration_section}. 

\subsection{The sets $H'$ and $H''$}

When we come out of the iteration loops (`reassignments') 
from the previous subsection, we finish with a set $S$ having 
a number of highly useful properties, among them:
\bigskip

$\bullet$ $|S|\ \geq\ |A|^{k - \delta}$;
\bigskip

$\bullet$ Each $R_y \subseteq R_x$; and,
\bigskip

$\bullet$ If we let $H$ denote those strings $h$ of length $k/2$ such that
$$
|R_h|\ \geq\ |A|^{k/2 - 2\delta},
$$
then for every such $h$ we will have that
$$
|\Sigma(S)|^{1-\eps/400c}\ <\ |\Sigma(R_h)|\ \leq\ |\Sigma(S)|.
$$
One can easily show, using the lower bound for $|S|$, 
that for $|A|$ sufficiently large,
$$
|H|\ >\ |A|^{k/2 - 2 \delta}.
$$
\bigskip

Since
$$
\sum_{z \in R_x} |\{h \in H\ :\ hz \in S \}|\ \geq\ 
|H|\cdot |A|^{k/2 - 2\delta},
$$
we deduce that there exists $z\in R_x$ such that there 
are at least 
$$
|H|\cdot |A|^{-2\delta}\ \geq\ |A|^{k/2 - 4\delta}  
$$
vectors $h \in H$ satisfying 
\begin{equation} \label{v_hold}
hz\ \in\ S.
\end{equation}
Fix one of these $z$, and let 
$$
H'\ \subseteq\ H
$$
denote all those $h \in H$ such that (\ref{v_hold}) holds.  Note that 
$$
|H'|\ \geq\ |A|^{k/2 - 4\delta}.
$$   

Next, let 
$$
H''\ \subseteq\ H'
$$
denote those $h \in H'$ such that there are at least 
\begin{equation} \label{H''count}
|H'|\cdot |\Sigma(H')|^{-1}/2 
\end{equation}
other $h' \in H'$ satisfying 
$$
\Sigma(h')\ =\ \Sigma(h).
$$
We have that 
$$
|H' \setminus H''|\ \leq\ |\Sigma(H')| (|H'| \cdot |\Sigma(H')|^{-1}/2)\ =\ 
|H'|/2
$$
So,
\begin{equation} \label{H''lower}
|H''|\ \geq\ |H'|/2\ \geq\ |A|^{k/2 - 5 \delta},
\end{equation}
for $|A|$ sufficiently large.

We also note that
$$
|\Sigma(H'')|\ \leq\ |\Sigma(H')|\ =\ 
|\Sigma( \{ hz\ :\ h \in H'\})|\ \leq\ |\Sigma(S)|.
$$
This is one of the places where it was essential to have that 
$z \in R_h$ for all $h \in H'$.
\bigskip

Now suppose that, in fact, 
\begin{equation} \label{H'_check}
|\Sigma(H'')|\ \leq\ |\Sigma(S)|^{1-\eps/400c}.
\end{equation}
If so, then we assign
$$
S\ \leftarrow\ H'',\ k\ \leftarrow\ k/2,\ \delta\ \leftarrow\ 5 \delta,
$$
and we repeat our iteration process again, starting in subsection
\ref{iteration_section}.
\bigskip

On the other hand, if (\ref{H'_check}) does not hold, then we will have
that 
\begin{equation} \label{H'H''}
|\Sigma(S)|^{1-\eps/400c}\ \leq\ |\Sigma(H'')|\ \leq\ 
|\Sigma(H')|\ \leq\ |\Sigma(S)|
\end{equation}

\subsection{The final leg of the proof}

From the fact that
$$
|\Sigma(\{hu\in S\ :\ h \in H'', u \in R_h \})|\ \leq\ |\Sigma(S)|,
$$
along with the fact that $R_h \subseteq R_x$ and
$$
|\Sigma(S)|^{1-\eps/400c}\ \leq\ |\Sigma(R_h)|\ \leq\ 
|\Sigma(R_x)|\ \leq\ |\Sigma(S)|,
$$
as well as (\ref{H'H''}), we deduce that there are at least 
$$
|\Sigma(S)|^{3 - 3\eps/400c}
$$
quadruples
$$
h_1, h_2\ \in\ \Sigma(H''),\ {\rm and\ } u_1, u_2\ \in\ \Sigma(R_x),
$$
such that
$$
\Sigma(h_1) + \Sigma(u_1)\ =\ \Sigma(h_2) + \Sigma(u_2).
$$

Now we apply Theorem \ref{BSG}, setting 
$$
X\ :=\ \Sigma(H''),\ {\rm and\ } Y\ :=\ \Sigma(R_x).
$$ 
Following the comment after Theorem \ref{SSV}, we have that there exists
$$
\Sigma\ \subseteq\ \Sigma(H''),\ |\Sigma|\ \geq\ |\Sigma(H'')|^{1-\eps/2c},
$$
such that 
\begin{equation} \label{small_double}
|\Sigma + \Sigma|\ \leq\ |\Sigma|^{1+\eps/2c}.
\end{equation} 
\bigskip

Let $H'''$ denote the set of all
$$
h\ \in\ H'',
$$
such that 
$$
\Sigma(h)\ \in\ \Sigma. 
$$
By (\ref{H''count}) and (\ref{H'H''}), we have that 
\begin{eqnarray*}
|H'''|\ &\geq&\ |\Sigma| (|H'|\cdot |\Sigma(H')|^{-1}/2)\\
&\geq&\ |\Sigma(H'')|^{1-\eps/2c} |H'|\cdot |\Sigma(S)|^{-1}/2\\
&\geq&\ |\Sigma(H'')|^{1-\eps/2c} |\Sigma(H'')|^{-1/(1-\eps/400c)} |H'|/2 \\
&\geq&\ |\Sigma(H'')|^{-\eps/c} |H'| \\
&\geq&\ |A|^{k/2 - 4 \delta - \eps}.
\end{eqnarray*}

By simple averaging, there is some vector 
$$
w\ \in\ A^{k/2-1},
$$
such that there are at least 
$$
|A|^{1 - 4 \delta - \eps}
$$
vectors $h \in H'''$ whose last $k/2-1$ coordinates are the vector $w$.
The upshot of this is that if we let 
$$
A'\ :=\ \{ a \in A\ :\ aw \in H'''\},
$$
then
\begin{equation} \label{A'_lower_bound}
|A'|\ \geq\ |A|^{1- 4 \delta - \eps},
\end{equation}
and
$$
A' + A' + 2 \Sigma(w)\ \subseteq\ \Sigma(H''') + \Sigma(H''')\ =\ 
\Sigma + \Sigma.
$$
\bigskip

Now we apply a weak form of the Ruzsa-Plunnecke Theorem
\cite{ruzsa}, given as follows:

\begin{theorem} \label{plunnecke_ruzsa} 
Suppose that $X$ is some finite subset of an additive abelian group, such that 
$$
|X + X|\ \leq\ C |X|.
$$
Then, we have that 
$$
|kX|\ =\ |X + X + \cdots + X|\ \leq\ C^k |X|.
$$
\end{theorem}

Using 
$$
X\ :=\ \Sigma,\ {\rm and\ } C\ :=\ |\Sigma|^{\eps/2c},
$$
we deduce that for $\ell$ even,
$$
|\ell A'|\ \leq\ |\ell \Sigma|\ \leq\ 
|\Sigma|^{1 + \eps \ell/2c}\ \leq\ |A|^{c + \eps \ell}
\ \leq\ |A'|^{(c + \eps \ell)/(1 - 4\delta -\eps)}
$$
By selecting $\delta > 0$ small enough, relative to $\eps > 0$,
we can ensure that for $\eps < 1/2$, 
$$
|\ell A'|\ \leq\ |A'|^{c(1 + 2 \eps \ell)}.
$$
Of course, when $1/2 \leq \eps < 1$ the inequality is trivial, as $c > 1$.
Clearly, on rescaling $\eps$ appropriately, our theorem is proved.

\section{Acknowledgements}

We would first and foremost like to thank Jozsef Solymosi, who emailed
us the reference to the Sudakov-Szemer\'edi-Vu paper, and then who
presented to us his own proof of it.  We would also like to thank 
Boris Bukh, who also pointed out this paper and result to us.  

Jozsef (as we have been reminded) and especially Boris made certain very 
helpful remarks to us concerning an application of the 
Bourgain-Chang sum-product theorem on a paper of ours yet to be 
written up, of which the present paper was once a part.  We wish to 
thank them for their remarks in advance of the writeup of that new paper.

Van Vu also mentioned his hypergraph result to us in passing during the
conference, and we wish to thank him for the comment.

We wish to thank Terry Tao for pointing out a result from section 2.6
from his book with Vu.  Evan plans to work on this.

We wish to thank Harald Helfgott for pointing out section 2.6 of the 
Tao-Vu book.  

Finally, we wish to thank Antal Balog for the warmth and interest in our
result.

For anyone we forgot to thank, our apologies in advance.

\end{document}